\documentclass[12pt]{article}
 \pdfoutput=1 
\usepackage{times,amssymb,amsmath,algorithm,algorithmicx,algpseudocode}
\usepackage{graphics,graphicx,psfrag}
\usepackage{amsthm}

\theoremstyle{plain}
\newtheorem{theorem}{Theorem}
\newtheorem*{problem}{Problem}

\theoremstyle{definition}
\newtheorem*{acknowledgement}{Acknowledgement}
\newtheorem{rem}[theorem]{Remark}

\newcommand{\h}{\mathcal{H}}

\def\eps{\varepsilon}

\title{The number of the maximal triangle-free graphs}
\author{J\'{o}zsef Balogh\thanks{ Department of Mathematics, University of Illinois, Urbana, IL 61801, USA and Bolyai Institute, University of Szeged, Szeged, Hungary {\tt jobal@math.uiuc.edu}.
 Research is partially supported by Simons Fellowship, NSF CAREER Grant DMS-0745185, Arnold O. Beckman Research Award (UIUC Campus Research Board 13039) and Marie Curie FP7-PEOPLE-2012-IIF 327763.} \and \v{S}\'{a}rka Pet\v{r}\'{i}\v{c}kov\'{a}\thanks{Department of Mathematics, University of Illinois, Urbana, IL 61801, USA {\tt petrckv2@illinois.edu}}}

\begin{document}
\maketitle

\begin{abstract}
Paul Erd\H{o}s suggested the following problem: Determine or estimate the number of maximal triangle-free graphs on $n$ vertices. Here we show that the number of maximal triangle-free graphs is at most $2^{n^2/8+o(n^2)}$, which matches the previously known lower bound. Our proof uses among others the Ruzsa-Szemer\'{e}di triangle removal lemma, and recent results on characterizing of the structure of independent sets in hypergraphs.

\end{abstract}

The maximum triangle-free graph has $n^2/4$ edges \cite{Mantel}. Hence, the number of triangle-free graphs is at least $2^{n^2/4}$, which was shown to be the correct order of magnitude by Erd\H{o}s, Kleitman and Rothschild~\cite{Erdos}  (see Balogh-Morris-Samotij~\cite{Balogh} or Saxton-Thomason~\cite{Saxton} for recent proofs). Moreover, almost every triangle-free graph is bipartite~\cite{Erdos}, even if there is a restriction on the number of edges (first shown by Osthus-Pr\"{o}mel-Taraz~\cite{Osthus}, extended by Balogh-Morris-Samotij-Warnke~\cite{Balogh2}; see ~\cite{Balogh} and ~\cite{Balogh2} for a more detailed history of the problem). This suggests that most of those graphs are bipartite, and subgraphs of a complete bipartite graph, therefore most of them are not maximal.  Erd\H{o}s suggested the following problem (as stated in Simonovits~\cite{Simonovits}): 
\begin{problem}
Determine or estimate the number of maximal triangle-free graphs on $n$ vertices.
\end{problem}

The following folklore construction (see \cite{Zare}, but it was known much earlier) shows that $f(n)\geq 2^{n^2/8+o(n^2)}$. Let $H$ be a graph on a vertex set $X\cup Y$ with $|X|=|Y|=n/2$ such that $X$ induces a perfect matching, $Y$ is an independent set, and there are no edges between $X$ and $Y$. For each pair of a matching edge $x_1 x_2\in E(H[X])$ and a vertex $y\in Y$, we add one of the edges $x_1 y$ or $x_2 y$ to $H$. Since there are $n/4$ matching edges in $E(H[X])$ and $n/2$ vertices in $Y$, we obtain $2^{n^2/8}$ triangle-free graphs, most of which are maximal.

In this note we prove a matching upper bound.

\begin{theorem} \label{thm_main}
The number of maximal triangle-free graphs with vertex set $[n]$ is at most $2^{n^2/8+o(n^2)}$.
\end{theorem}

Our first tool is  a corollary of recent powerful counting theorems of Balogh-Morris-Samotij~\cite[Theorem 2.2.]{Balogh}, and Saxton-Thomason~\cite{Saxton}.

\begin{theorem}\label{thm1} 
For each $\delta>0$ there is $t<2^{O(\log n\cdot n^{3/2})}$ and a set $\{G_1,\dots,G_t\}$ of graphs, each containing   at most $\delta n^3$ triangles, such that for every triangle-free graph $H$ there is $i\in [t]$ such that $ H \subseteq  G_i$, where   $n$ is sufficiently large.
\end{theorem}

By the Erd\H{o}s-Simonovits supersaturation theorem~\cite{Erdos2}, each $G_i$ has at most about $n^2/4$ edges.

\begin{theorem}\label{thm_saturation} For every $\gamma >0$ there is $\delta(\gamma)>0$ such that every $n$-vertex graph with $(1/4+\gamma)n^2$ edges contains as least $\delta(\gamma) n^3$ triangles.
\end{theorem}

We also use the Ruzsa-Szemer\'{e}di triangle-removal lemma~\cite{Ruzsa}.
\begin{theorem}\label{thm2}
For every $\eps>0$ there is $\delta(\eps)>0$ such that any graph $G$ on $n$ vertices with at most $\delta(\eps) n^3$ triangles can be made triangle-free by removing at most $\eps n^2$ edges.
\end{theorem}

Our next tool is the following theorem of Hujter and Tuza~\cite{Hujter}. Recall that a set $I\in V(G)$ is an \emph{independent set} if no two vertices in $I$ are adjacent. An independent set $I$ is a \emph{maximal independent set} if $I\cup \{v\}$ contains an edge for every $v\in V(G)- I$. Note that we write $|G|$ for the number of vertices of $G$.

\begin{theorem}\label{thm3}
Every triangle-free graph $G$ has at most $2^{|G|/2}$ maximal independent sets. 
\end{theorem}

In the next section we prove our main result, Theorem~\ref{thm_main}.

\section[Proof of Main Theorem]{Proof of Theorem~\ref{thm_main}}

We show that for every $\eps>0$ and for every $\gamma>0$, the number of maximal triangle-free graphs with vertex set $[n]$ is $2^{(1/8+2\eps+\gamma)n^2}$ for sufficiently large $n$. We fix  arbitrarily small constants $\eps, \gamma>0$. First we  apply Theorem~\ref{thm2} with this $\eps$, and Theorem~\ref{thm_saturation} with this $\gamma$, which provides us $\delta(\eps)$ and $\delta(\gamma)$.
We define $\delta=\min\{\delta(\eps), \delta(\gamma)\}$ and   then apply Theorem~\ref{thm1} with this $\delta$. For every $i\in[t]$, we count the number of maximal triangle-free graphs $H$ that satisfy $ H\subseteq  G_i$. Denote $\h$  the set of  maximal triangle-free graphs with vertex set $[n]$, and let $\h_i=\{H\in\h:  H\subseteq G_i\}$. 

Since $t \leq 2^{\eps n^2}$ for sufficiently large $n$, we have 
$$ |\h|\leq \sum_{i=1}^t  |\h_i| \leq  2^{\eps n^2}\max_{i\in[t]}|\h_i|.$$
Fix an arbitrary $i\in [t]$.
 By Theorem~\ref{thm2} applied on $G_i$, there is $F_i\subseteq E(G_i)$ such that $|F_i|\leq \eps n^2$ and $G_i-F_i$ is triangle-free. For each $G_i$ we fix one such $F_i$. For every $F^*\subseteq F_i$ define  $\h_{i}(F^*)=\{H\in \h_i: E(H)\cap F_i=F^*\}$.

\begin{figure}
\begin{center}
\includegraphics[width=\textwidth]{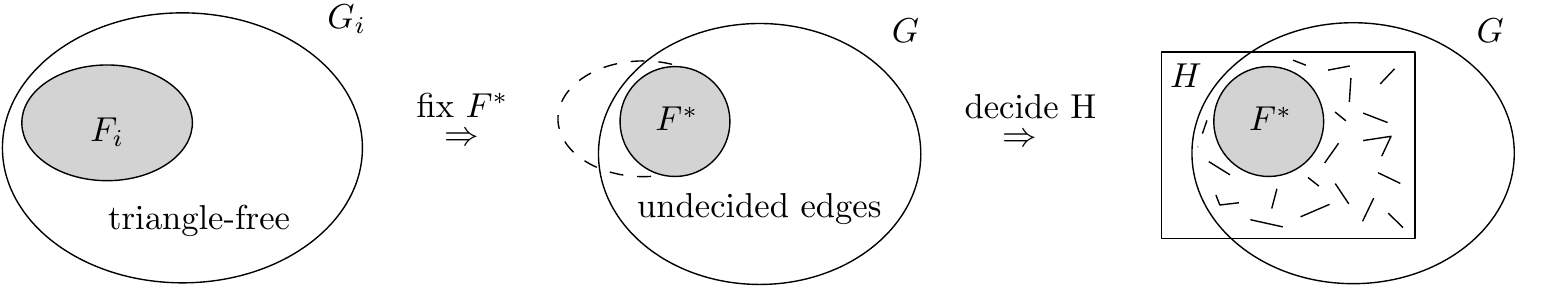}
\end{center}
\caption{The structure of $G_i$ and $F_i$.}
\end{figure} 
 
 Now we show that for every choice of $F^*$ we have $|\h_{i}(F^*)|\leq 2^{e(G_i)/2}$. Fix $F^*$, and let
   $$G:=G_i-(F_i-F^*)- \{e\in E(G_i): \exists f,g \in  F^* \mbox{ such that } e,f,g \mbox{ form a triangle} \}.$$

   So, $G$ is obtained from $G_i$ by removing edges that are in none of $H\in \h_{i}(F^*)$.  We can assume that $F^*$ is triangle-free since otherwise $\h_{i}(F^*)=\emptyset$. We now count the number of ways to add edges of $E(G)-F^*$ to $F^*$ such that the resulting graph is maximal triangle-free. We construct an auxiliary graph $T$ as follows:
$$
\begin{array}{lllll}
V(T):=E(G)-F^* \quad \mbox{ and }\quad
E(T):=\{ef| \mbox{ $\exists d\in F^*$: } \{d,e,f\} \mbox{ spans a triangle in } G \}.
\end{array}
$$

\noindent
{\bf Claim 1.} $T$ is triangle-free. 
\begin{proof}
Suppose not. Let $e,f,g$ be vertices of a triangle in $T$.  Then $e,f,g\in E(G)- F^*$ and there are $d_1, d_2, d_3\in F^*$ such that the $3$-sets $\{d_1,e,f\}$, $\{d_2,e,g\}$, and $\{d_3,f,g\}$ span triangles in $G$. 
As $G_i - F_i$ is triangle-free and $G - F^* \subseteq G_i - F_i$, it follows that the edges $e,f,g$ share a common endpoint in $G$, and that  $\{d_1,d_2,d_3\}$ spans a triangle. This is a contradiction since $F^*$ is triangle-free. 
\end{proof}

\noindent
{\bf Claim 2.}  If $H\in \h_{i}(F^*)$, then $E(H)- F^*$ spans a maximal independent set in $T$.

\begin{proof} Let $H\in \h_{i}(F^*)$. We first show that $E(H)-F^*$ spans an independent set in $T$. If not, then there is an edge $ef$ in $E(T)$ with $e,f\in E(H)-F^*$. By the definition of $E(T)$, there is $d\in F^*$ such that the edges $d,e,f$ form a triangle in $G$, which is clearly in $H$.

Suppose now that $E(H)- F^*$ is an independent set in $T$ that is not maximal. So, there is $x\in E(G)- E(H)$ such that for every $y\in E(H)-F^*$ and for every $z\in F^*$, the edges $x,y,z$ do not span a triangle in $G$. This means that $H\cup \{x\}$ is triangle-free. Hence, $H$ is not maximal.
\end{proof}

By Theorem~\ref{thm3}, the number of maximal independent sets in $T$ is at most $2^{|T|/2}$. Thus

$$|\h_{i}(F^*)|\leq 2^{|T|/2}\leq 2^{e(G_i)/2}\leq 2^{(n^2/4+\gamma)/2},$$
where the last inequality follows from Theorem~\ref{thm_saturation}.

The number of ways to choose $F^*\subseteq F_i$ for a  given $F_i$ is at most $2^{|F_i|}\leq 2^{\eps n^2}$, so we can conclude that for sufficiently large $n$,

$$ |\h|\leq  2^{\eps n^2}\max_{i\in[t]}|\h_i|\leq  2^{\eps n^2}\sum_{F^*\subseteq F_i} |\h_{i}(F^*)|\leq 2^{\eps n^2}2^{\eps n^2} \max_{F^*\subseteq F} |\h_{i}(F^*)|$$
$$\leq 2^{2\eps n^2} 2^{(n^2/4+\gamma n^2)/2} \leq  2^{(1/8 +2\eps+\gamma)n^2}.$$

\section{Concluding remarks}

It would be interesting to have similar results for $K_{r+1}$ as well. A straightforward generalization of the construction for maximal triangle-free graphs implies that there are at least $r^{n^2/(4r^2)}$ maximal $K_{r+1}$-free graphs. Unfortunately, not all steps of our upper bound method work when $r>2$. We are able to get only the following modest improvement on the trivial $2^{(1-1/r+o(1))n^2/2}$ bound: for every $r$ there is a positive constant  $c_r$ such that the number of maximal $K_{r+1}$-free graphs is at most $2^{(1-1/r-c_r)n^2/2}$ for $n$ sufficiently large. More precisely, if we let $s=2^{\binom{r+1}{2}-1}$, then the number of maximal $K_{r+1}$-free graphs is at most $(s-1)^{n^2/(r(r+2))+o(n^2)}$.

A similar question was raised by Cameron and Erd\H{o}s in \cite{CE}, where they asked how many maximal sum-free sets are contained in $[n]$. They were able to construct $2^{n/4}$ such sets. An upper bound $2^{3n/8+o(n)}$ was proved by Wolfovitz~\cite{wolfovitz}. Our proof method instantly improves this upper bound to $3^{n/6+o(n)}$, as observed in~\cite{team}. Balogh-Liu-Sharifzadeh-Treglown~\cite{team} pushed the method further to prove a matching upper bound, $2^{n/4+o(n)}$. As~\cite{team}  contains all the details, we omit further discussion here.

\section*{Recent development}    
\begin{rem}
 Alon pointed out that if the number of $K_r$-free graphs is $2^{c_r n^2 +o(n^2)}$, then $c_r$ is monotone (though not clear if strictly monotone) increasing in $r$.
\end{rem}  

\begin{rem} 
A discussion with Alon and \L uczak at IMA led to the following construction that gives $2^{(1-1/r)n^2/4+o(n^2)}$ maximal $K_{r+1}$-free graphs: partition the vertex set $[n]$ into $r$ equal classes, place a perfect matching into $r-1$ of them. Between the classes we have the following connection rule: between two matching edges place exactly $3$ edges, and between a vertex (from the class which is an independent set) and a matching edge put exactly $1$ edge.

\end{rem}  

\begin{rem}
In yet to be published work, Balogh, Liu, Pet\v{r}\'{i}\v{c}kov\'{a}, and Sharifzadeh proved that almost every maximal
triangle-free graph $G$ admits a vertex partition $X\cup Y$ such that $G[X]$ is a perfect matching and $Y$ is an independent set, as in the construction.
\end{rem}

\begin{acknowledgement} 
We thank Rupei Xu for showing us the construction for the lower bound. We thank the referee for a quick review and pointing out some typoes.
\end{acknowledgement}

\end{document}